\documentclass[final,1p,times]{elsarticle}

\usepackage{amssymb}
\usepackage{float}
\usepackage{balance}

\usepackage{amsthm}
\usepackage{amsmath}
\usepackage{natbib}
\usepackage{epstopdf}
\usepackage[colorlinks,
            linkcolor=red,
            anchorcolor=blue,
            citecolor=green
            ]{hyperref}
\usepackage{multirow}

\newtheorem{lemma}{Lemma}
\newtheorem{remark}{Remark}
\newtheorem{theorem}{Theorem}

\begin{document}

\begin{frontmatter}

\title{ A novel multi-step method for the partial pole assignment in symmetric quadratic pencil with time delay}

\author[mymainaddress,mysecondaryaddress]{Qing Liu\corref{mycorrespondingauthor}}
\cortext[mycorrespondingauthor]{Corresponding author.}
\ead{qliu@seu.edu.cn}

\address[mymainaddress]{School of Mathematics and Shing-Tung Yau Center, Southeast University, Nanjing 211189, P. R. China.}
\address[mysecondaryaddress]{Nanjing Center for Applied Mathematics, Nanjing 211135, P. R. China.}

\begin{abstract}
In this paper, we study the partial pole assignment problem in symmetric quadratic pencil with time delay.
A novel multi-step method is proposed to solve this problem,  resulting in the undesired eigenvalues being moved to desired values, and the remaining eigenvalues  unchanged.
By establishing a new matrix equality relation and using a multi-step method,  the problem  is transformed  into solving linear systems with low order.  Specifically, assuming that there are $p$ undesired eigenvalues requiring reassigned,  the size of the linear system we finally solved is $p^2.$
Notably, the method demonstrates high efficiency for large systems with only a few poles requiring reassigned.
Numerical examples are provided to illustrate the effectiveness of the proposed method.\\
\end{abstract}

\begin{keyword}
Inverse problem; Second-order linear system; Partial pole assignment; Multi-step method; Time delay.
\end{keyword}

\end{frontmatter}

\section{ Introduction}

\renewcommand{\baselinestretch}{1.26}\large\normalsize
In structural dynamics (e.g., bridges, buildings, airplanes and automobiles), the finite element method is used to discretize the damped linear system with $n$ degrees of freedom, and the second-order linear differential equation with constant coefficients is generated as follows:
\begin{equation}\label{eq1.1}
M{\kern 1pt} \ddot x(t) + C{\kern 1pt} \dot x(t) + Kx(t) = f(t-\tau),
\end{equation}
where $M,C,K \in {\mathbb{R}^{n \times n}}$  are the symmetric mass matrix, damping matrix and stiffness matrix with $M$ positive definite, $x(t) \in \mathbb{R}^{n}$ is the displacement vector, $f(t-\tau)$ is the external force control or control vector and $\tau$ is the time delay between the measurements of the state and the actuation of the control. Let $f(t-\tau)=0$, then \eqref{eq1.1} may be solved by using the method of separation of variables. Substituting $x(t) = x{e^{\lambda t}}$, where $\lambda  \in {\mathbb{C}}$, $x \in {\mathbb{C}^n}$, into \eqref{eq1.1} yields the quadratic eigenvalue problem of the open-loop system:
\begin{equation}\label{eq1.2}
P(\lambda )x: = \left ({\lambda ^2}M + \lambda C + K\right )x = 0,
\end{equation}
where $P(\lambda ) = {\lambda ^2}M + \lambda C + K$ is called the open-loop quadratic pencil,  $\lambda$ is referred to as the eigenvalue or complex frequency or pole, and  $x$ is called the eigenvector or vibration mode. In practical engineering problems, only a few poles are unstable and lead to the system resonance. In order to prevent system instability or resonance, it is essential to modify the unstable poles of $P(\lambda ) $ through feedback control. At the same time, it is crucial to ensure that the rest of the system's eigenstructure remains unaltered, thereby preserving the property of no spill-over. This leads to the partial quadratic pole assignment problem \cite{1,2,3,4}.

 The system can be controlled by the external force $f(t-\tau) = Bu(t-\tau)$, and the corresponding equation of state is described by
\begin{equation}\label{eq1.3}
M{\kern 1pt} \ddot x(t) + C{\kern 1pt} \dot x(t) + Kx(t) =Bu(t-\tau),
\end{equation}
where $B \in {\mathbb{R}^{n \times m}}$ is the full column rank control matrix and $u(t-\tau)\in {\mathbb{R}^{ m}}$  is the feedback control vector. If $m =1$, the system \eqref{eq1.3} is reduced into the single-input control system, and if $m > 1$, it is called the multi-input control system. In this paper, $u(t-\tau)$ is chosen as the displacement-velocity feedback control
\begin{equation}\label{eq1.4}
u(t-\tau) = G^\top x(t-\tau)+{F^\top}\dot x(t-\tau),
\end{equation}
where $G, F \in \mathbb{R}^{n\times m}$ are the displacement feedback matrix and velocity feedback matrix, respectively. Substituting \eqref{eq1.4} into \eqref{eq1.3},   the following second-order closed-loop system can be obtained:
\begin{equation}\label{eq1.5}
M{\kern 1pt} \ddot x(t) +  C{\kern 1pt}  \dot x(t)  + K{\kern 1pt} x(t)  = B\left({F^\top}\dot x(t-\tau)  + {G^\top}x(t-\tau) \right) .
\end{equation}
It follows that the corresponding quadratic eigenvalue problem with time delay is
\begin{equation}\label{eq1.6}
{P_\tau}(\lambda )x: = \left [{\lambda ^2}M + \lambda \left (C - B{F^\top}e^{-\lambda \tau} \right) + \left( K - B{G^\top}e^{-\lambda \tau}\right)\right ]x = 0,
\end{equation}
where ${P_\tau}(\lambda )= {\lambda ^2}M + \lambda \left (C - B{F^\top}e^{-\lambda \tau} \right) + \left( K - B{G^\top}e^{-\lambda \tau}\right)$.
Suppose that $\left\{ {{\lambda _i}} \right\}_{i = 1}^p(p< 2n)$ are undesired poles of the second-order system \eqref{eq1.3}, which need to be reassigned to the desired poles $\left\{ {{\mu _i}} \right\}_{i = 1}^p$. This leads to the following problem, known as the partial quadratic pole assignment problem with time delay.

\vspace{0.2cm}

\textbf{Problem 1.} Given the symmetric matrices of the second-order system $M,C,K\in \mathbb{R}^{n\times n}$ with $M$ positive definite, the full column rank control matrix $B \in {\mathbb{R}^{n \times m}}$, the time delay $\tau> 0 $, the self-conjugate subset $\left\{ {{\lambda _i}} \right\}_{i = 1}^p\left( {p < 2n} \right)$ of the open-loop spectrum $\left\{ {{\lambda _i}} \right\}_{i = 1}^{2n}$ and the corresponding eigenvector set $\left\{ {{x_i}} \right\}_{i = 1}^p$, and the desired self-conjugate set $\left\{ {{\mu_i}} \right\}_{i = 1}^p$. Find the state feedback matrices $F,G$ $\in {\mathbb{R}^{n \times m}}$ such that the closed-loop delayed pencil ${P_\tau}(\lambda )$ has the desired eigenvalues $\left\{ {{\mu_i}} \right\}_{i = 1}^p$ and the eigenpairs $\left\{ {\left( {{\lambda _i},{x_i}} \right)} \right\}_{i = p + 1}^{2n}$.

If the system \eqref{eq1.3} is partially controllable \cite{5} with respect to the eigenvalues $\lambda=\lambda_i (i=1,2,\ldots,p)$, then
$${\rm{rank}}\big\{ {{\lambda ^2}M + \lambda C + K,B} \big\} = n.$$
In this paper, we assume that $\left\{ {{\lambda _i}} \right\}_{i = 1}^p$ and $\left\{ {{\mu _i}} \right\}_{i = 1}^p$ are all distinct, the system \eqref{eq1.3} is partially controllable with respect to $\left\{ {{\lambda _i}} \right\}_{i = 1}^p$, and
\begin{equation}\label{eq1.7}
\left\{\mu_i \right\}_{i=1}^p \cap \left\{\lambda_i \right\}_{i=1}^{2n} = \emptyset, \
\left\{\lambda_i \right\}_{i = 1}^p \cap \left\{\lambda_i \right\}_{i=p+1}^{2n} = \emptyset.
\end{equation}
Let
$$\Lambda  ={\rm {diag}}\big(\lambda_1,\lambda_2, \ldots ,\lambda_{2n}\big),\  X = \big[x_1,x_2, \ldots ,x_{2n}\big],$$
$$\Lambda_1 = {\rm {diag}}\big(\lambda_1,\lambda_2, \ldots, \lambda_p \big), \ \Lambda_2 = {\rm {diag}}\big(\lambda_{p+1},\lambda_{p+2}, \ldots,\lambda_{2n}\big),$$
$$X_1 = \big[x_1, x_2, \ldots , x_p\big], \  X_2 = \big[x_{p+1},x_{p+2}, \ldots ,x_{2n}\big].$$
 Datta et al \cite{6} derived the following three orthogonality relations
\begin{equation}\label{eq1.8}
\Lambda X^\top MX \Lambda -X^\top KX=D_1,
\end{equation}
\begin{equation}\label{eq1.9}
\Lambda X^\top CX \Lambda + \Lambda X^\top KX + X^\top KX \Lambda =D_2,
\end{equation}
\begin{equation}\label{eq1.10}
\Lambda X^\top MX + X^\top MX \Lambda +X^\top CX =D_3,
\end{equation}
where $D_1,D_2,D_3 \in {\mathbb{C}^{2n \times 2n}}$ are diagonal matrices. Via orthogonality relation \eqref{eq1.8}, Datta et al \cite{6} presented the explicit solution of feedback control matrices $F$ and $G$ for the partial quadratic pole assignment problem for the single-input system. Ram and Elhay \cite{7} generalized it to the multi-input system such that a multi-step method for the partial quadratic pole assignment problem is developed. However, this method cannot reassign zero or near zero poles and needs to solve $(m\cdot p)$ $n$-order linear systems. Based on the receptance method \cite{8}, Ram et al \cite{9} proposed an  algorithm for the partial quadratic pole assignment problem with time delay for the single-input system. Bai et al \cite{10,11}, Ram and Mottershead \cite{12}, Chen and Xie \cite{18}, and Richiedei et al \cite{20} developed the receptance methods for the partial quadratic pole assignment problems for the multi-input system. However, the receptance matrix needs to be measured by modal tests, and the receptance matrix is sensitive to the test noise \cite{8,13}. Mao \cite{14} provided a method for solving the partial quadratic pole assignment problem without using the receptance matrix. However, this method needs to solve the Sylvester equation. Liu and Yuan \cite{15} developed a multi-step method for solving the partial quadratic pole assignment problem with time delay. However, this method needs to solve $(m\cdot p)$ $n$-order linear systems, which the computational cost is relatively expensive for the large-scale system.

 In this paper, we study an efficient numerical method for the partial pole assignment in symmetric quadratic pencil with time delay.  The main innovations presented in this paper can be summarized as follows.  Firstly, we provide an explicit solution to \textbf{Problem 1} for the single-input system, utilizing the third orthogonality relation \eqref{eq1.10} and the inverse of the Cauchy matrix. This method involves only the multiplication of matrices and vectors. Secondly, for the multi-input system, we establish a new connection between system matrices, undesired eigenvalues and corresponding eigenvectors, as well as desired eigenvalues and control vectors. This enables us to develop a novel multi-step method that transforms \textbf{Problem 1} into solving $p^{2}$-order linear systems,  and does not involve solving  $n$-order linear systems.  This approach proves to be highly efficient for large-scale systems, as it only requires reassigning a small number of poles, i.e. it is particularly suitable for $p\ll n$.

% by establishing the new connection among system matrices, undesired eigenvalues and corresponding eigenvectors, desired eigenvalues and control vectors, we develop a novel multi-step method to transform this problem into solving $m-1$ $p^{2}$-order linear systems, where $p\ll n$, which is  quite efficient for large systems where only a small number of poles need to be reassigned.

The remaining of this paper is organized as follows. In section \ref{single}, we present an explicit solution to \textbf{Problem 1} for the single-input system. In section \ref{multi},
a  novel multi-step method is proposed to solve \textbf{Problem 1} for the  multi-input system.  Numerical examples are presented in
section \ref{example}. Finally, we draw some conclusions  in section \ref{con}.

\section{Single-input control }\label{single}
\hspace{2pt}
In this section, we consider the single-input control system, which means $m=1$,
$$
B = b \in {\mathbb{R}^n}, \ F = f \in {\mathbb{R}^n}, \ G = g \in {\mathbb{R}^n}.
$$
To begin with, we  provide the following three lemmas.

 \begin{lemma}\label{lemma1.1}\cite{5}
 {The single-input control system \eqref{eq1.3} is  partially controllable with respect to the eigenvalues $ {{\lambda _i}}$ if and only if  $b^\top x_i \neq 0 \ $ $(i=1,2, \ldots, p)$.}
 \end{lemma}

 According to the orthogonality relation \eqref{eq1.10}, we can obtain
\begin{equation}\label{eq2.11}
\Lambda_1 X_1^\top MX_2 + X_1^\top MX_2 \Lambda_2 +X_1^\top CX_2 =0.
\end{equation}
Therefore, the following lemma can be verified.

 \begin{lemma}\label{lemma1.2}\cite{10}
 For arbitrary $\beta \in \mathbb{C}^p$, define
\begin{equation}\label{eq2.12}
f = M{X_1}\beta, \   g = \left(M{X_1}{\Lambda_1} + C{X_1}\right)\beta ,
\end{equation}
then the closed-loop delayed pencil ${P_\tau}(\lambda )$ has the eigenpairs $\left\{ {\left( {{\lambda _i},{x_i}} \right)} \right\}_{i = p + 1}^{2n}$, that means
$$
M{X_2}\Lambda_2^2 + C{X_2}{\Lambda_2}- bf^\top{X_2}{\Lambda_2}e^{-\tau\Lambda_{2} } + K{X_2} - bg^\top{X_2}e^{-\tau\Lambda_{2} } = 0.
$$
 \end{lemma}

 \begin{remark}\label{remark1}
 Lemma \ref{lemma1.2} indicates that the feedback vectors $f$ and $g$ defined by \eqref{eq2.12} can make the closed-loop delayed pencil ${P_\tau}(\lambda )= {\lambda ^2}M + \lambda \left(C - b{f^\top}e^{-\lambda \tau} \big ) + \big( K - b{g^\top}e^{-\lambda \tau}\right)$ keep the remaining eigenstructure of the original system unchanged, i.e., the no spill-over property is preserved.
 \end{remark}

 Assume that $\left\{(\mu_i,y_i) \right\}_{i=1}^p$ are the eigenpairs of the closed-loop delayed pencil ${P_\tau}(\lambda )$, define
$$
\Sigma_1 ={\rm{ diag}} \big(\mu_1, \mu_2, \ldots, \mu_p\big), \ Y_1 = \big[y_1, y_2, \ldots, y_p\big],
$$
\begin{equation}\label{eq2.13}
H = \big(h_{ij} \big) =\Lambda_1X_1^\top MY_1+X_1^\top MY_1\Sigma_1+X_1^\top CY_1.
\end{equation}

 \begin{lemma}\label{lemma1.3}\cite{15} If the single-input control system \eqref{eq1.3} is  partially controllable with respect to the eigenvalues $\left\{\lambda_i \right\}_{i=1}^p$, $\left\{ {{\lambda _i}} \right\}_{i = 1}^p$ and $\left\{ {{\mu _i}} \right\}_{i = 1}^p$ are all distinct and satisfy \eqref{eq1.7}, and $\left\{y_i \right\}_{i=1}^p$ are the solutions of the linear equations
$$
\left(\mu_i^2M+\mu_i C+K\right )y_i=b, \ i=1,2, \ldots, p ,
$$
then $H$ is nonsingular, and
$$
H = \tilde H \hat H,
$$
where
\begin{equation}\label{eq2.98}
\tilde H = {\rm{diag}}(b^\top x_1, {b^\top}x_2, \ldots, b^\top x_p), \
\hat H = \left( {\begin{array}{*{20}{c}}
\frac{1}{\mu_1 - \lambda_1} & \cdots & \frac{1}{\mu_p - \lambda_1} \\
    \cdots                  & \cdots & \cdots                      \\
\frac{1}{\mu_1 - \lambda_p} & \cdots & \frac{1}{\mu_p - \lambda_p}
\end{array}} \right).
\end{equation}
 \end{lemma}

From Lemma \ref{lemma1.2}, solving \textbf{Problem 1} is transformed into determining vector $\beta$ in \eqref{eq2.12} such that
\begin{equation}\label{eq2.99}
MY_1\Sigma_1^2 + CY_1\Sigma_1 - bf^\top Y_1\Sigma_1e^{-\tau\Sigma_{1} }  + KY_1    - bg^\top Y_1e^{-\tau\Sigma_{1} }  = 0.
\end{equation}
Substituting  $f$ and $g$ defined in \eqref{eq2.12} into \eqref{eq2.99}, and according to the definition of $H$ in \eqref{eq2.13}, we can derive
\begin{displaymath}
\begin{aligned}
MY_1\Sigma_1^2 + CY_1\Sigma_1 + KY_1 &= b\beta^\top \left(\Lambda_1X_1^\top MY_1 + X_1^\top MY_1\Sigma_1 + X_1^\top CY_1\right)e^{-\tau\Sigma_{1} }\\
&= b\beta^\top H e^{-\tau\Sigma_{1} }: = b\gamma ,
\end{aligned}
\end{displaymath}
where $\gamma  = \beta^\top H e^{-\tau\Sigma_{1} }$. Moreover, the vector $\gamma$ is related to the eigenvector $Y_{1}$ of the closed-loop delayed pencil ${P_\tau}(\lambda )$. In order to obtain $Y_{1}$, we can calculate the eigenvectors $\left\{y_i \right\}_{i=1}^p$ from the linear equations
\begin{equation}\label{eq2.14}
\left(\mu_i^2M+\mu_i C+K\right )y_i=b, \ i=1,2, \ldots, p,
\end{equation}
which corresponds to the selection $\gamma=[1,1,\ldots,1] \in \mathbb{R}^{p}$. Therefore, parameter vector $\beta$ satisfies
\begin{equation}\label{eq2.15}
\beta^\top H = [1,1,\ldots,1]e^{\tau\Sigma_{1} }.
\end{equation}
According to Lemma \ref{lemma1.3}, the linear system \eqref{eq2.15} has a unique solution, then we can obtain the following result.

\begin{theorem}\label{th1}
If the single-input control system \eqref{eq1.3} is partially controllable with respect to the eigenvalues $\left\{\lambda_i \right\}_{i=1}^p$, $\left\{ {{\lambda _i}} \right\}_{i = 1}^p$ and $\left\{ {{\mu _i}} \right\}_{i = 1}^p$ are all distinct and satisfy \eqref{eq1.7}, and $\left\{y_i \right\}_{i=1}^p$ are the solutions of \eqref{eq2.14}, then \textbf{Problem 1} has a unique solution, which can be expressed by \eqref{eq2.12}, where $\beta$ is the solution of the linear system \eqref{eq2.15}.
\end{theorem}

In fact, we  need not to solve $p$ linear systems \eqref{eq2.14}. Since $\hat H$ defined in \eqref{eq2.98} is a Cauchy matrix, define $T = \big(t_{ij}\big)_{p \times p} = \hat H^{-1}$ and $\beta  =\big [\beta_1, \beta_2, \ldots, \beta_p\big]^\top$, based on the inverse of the Cauchy matrix \cite{16}, we can derive
\begin{equation}\label{eq2.16}
t_{ij} = \frac{\prod\limits_{k=1}^p (\lambda_j - \mu_k) \prod\limits_{k=1}^p (\mu_i - \lambda_k)}
{(\lambda_j - \mu_i)\prod\limits_{k=1,k \ne i}^p (\mu_i - \mu_k)\prod\limits_{k=1,k \ne j}^p (\lambda_j - \lambda_k)}.
\end{equation}
Therefore, the explicit expression of $\beta$ in \eqref{eq2.12} is given by
\begin{equation}\label{eq2.97}
\beta_j = \frac{1}{b^\top x_j}\sum\limits_{i=1}^p t_{ij} e^{\tau\mu_{i}}.
\end{equation}
Substituting $t_{ij}$ defined in \eqref{eq2.16} into \eqref{eq2.97}, we can obtain
\begin{equation}\label{eq2.17}
\beta_j = \frac{1}{b^\top x_j}\frac{\prod\limits_{k=1}^p (\lambda_j - \mu_k)}
{\prod\limits_{k=1,k \ne j}^p (\lambda_j - \lambda_k)}\sum\limits_{i=1}^p \frac{\prod\limits_{k=1}^p (\mu_i - \lambda_k)}
{(\lambda_j - \mu_i)\prod\limits_{k=1,k \ne i}^p (\mu_i - \mu_k)}e^{\tau\mu_{i}}.
\end{equation}

 \begin{remark}\label{remark2}The explicit expression of $\beta$ in \eqref{eq2.12} is easy to implement without solving any linear system. If time delay $\tau = 0$, \eqref{eq2.17} can be further simplified. Datta et al \cite{6} prove that
$$
\sum\limits_{i=1}^p \frac{\prod\limits_{k=1,k \ne j}^p (\mu_i - \lambda_k)}{\prod\limits_{k=1,k \ne i}^p (\mu_i - \mu_k)} = 1,
$$
then,
$$
\beta_j = \frac{1}{b^Tx_j}(\mu_j - \lambda_j)\prod\limits_{k=1,k \ne j}^p \frac{\lambda_j - \mu_k}{\lambda_j - \lambda_k}.
$$
 \end{remark}

 It is easy to prove that if the complex eigenvalues of $\Lambda_1$ and $\Sigma_1$ appear as conjugate pairs, the velocity feedback vector $f$ and displacement feedback vector $g$ calculated by this method are both real.
Based on Theorem \ref{th1}, we can summarize the following algorithm.

\vspace{0.5cm}

\begin{tabular}{l}
\hline
{\bf Algorithm 1} A direct method for \textbf{Problem 1} by single-input control. \\
\hline
{\bf Input: } \\
\indent  The symmetric matrices $M,C\in \mathbb{R}^{n\times n}$ with $M$ positive definite;\\
\indent  The control vector $b \in {\mathbb{R}^{m}}$ and the time delay $\tau > 0$;\\
\indent  The self-conjugate eigenpairs $\left\{(\lambda_i,x_i)\right\}_{i=1}^p$ and a self-conjugate set $\left\{ {{\mu_i}} \right\}_{i = 1}^p$.\\

{\bf Output: } \\
\indent State feedback real vectors $f$ and $g$. \\
1: Form $\Lambda_1 = {\rm{diag}}(\lambda_1,\lambda_2, \ldots ,\lambda_p)$ and $X_1 = [x_1,x_2, \ldots ,x_p]$;\\
2: For $j = 1,2, \ldots, p$, calculate $\beta_j $ according to \eqref{eq2.17};\\
3: Form $\beta  = \big[\beta_1,\beta_2, \ldots ,\beta_p\big]^\top$;\\
4: Compute $ f = MX_1\beta, \ g = \left(MX_1\Lambda_1 + CX_1\right)\beta$. \\
\hline
\end{tabular}

\section{ Multi-input control }\label{multi}

In this section, we consider \textbf{Problem 1} for the multi-input control system, which means $m>1$. The closed-loop system \eqref{eq1.5} can be written as
$$
M \ddot x(t) + C \dot x(t) + Kx(t) = \sum\limits_{k=1}^m b_k \left(f_k^\top \dot x(t-\tau) + g_k^\top x(t-\tau)\right),
$$
where $b_k$, $f_k$, and $g_k$ are the $k$-th columns of $B$, $F$ and $G$, respectively. The corresponding closed-loop delayed pencil is given by
$$
P_{\tau }(\lambda) = \lambda^2M + \lambda \Big(C - \sum\limits_{k=1}^m b_kf_k^\top e^{-\lambda \tau}\Big) + K - \sum\limits_{k=1}^m b_kg_k^\top e^{-\lambda \tau}.
$$
For $j = 1, 2, \ldots ,p$ and $k = 1, 2, \ldots ,m$, define
\begin{equation}\label{eq3.18}
 \xi_{jk} = \eta_{jk} + \frac{k}{m}\big(\mu_j - \eta_{jk}\big), \ \xi_{jm} = \mu_j.
\end{equation}
Let
\begin{equation}\label{eq3.19}
\begin{aligned}
&C_1 = C, \ C_k = C - \sum\limits_{i=1}^{k-1} b_if_i^\top e^{-\lambda \tau}, \\
&K_1 = K, \ K_k = K - \sum\limits_{i=1}^{k-1} b_ig_i^\top e^{-\lambda \tau} ,
\end{aligned}
\end{equation}
where $k = 2,3, \ldots, m$, then \textbf{Problem 1}  can be transformed into  partial pole assignment problems for $m$ single-input second-order linear systems. For $k = 1$, $f_{1}$ and $g_{1}$ can be obtained by Algorithm 1 such that the eigenvalues $\left\{ {{\lambda _i}} \right\}_{i = 1}^{p}$ of the open-loop pencil $P(\lambda )$ are reassigned to $\big\{\xi_{j1}\big\}_{j=1 }^p$. However, the system matrices $C_{2}$ and $K_{2}$ are no longer symmetric after the first-step assignment. Therefore, Algorithm 1 cannot be applied for subsequent assignment. For multi-input control systems, we  propose the following multi-step method to solve \textbf{Problem 1}. Suppose that after  $(k-1)$-step assignment, where $2\leq k<m$, the system
\begin{equation}\label{eq3.20}
M \ddot x(t) + C_k \dot x(t) + K_k x(t) = 0
\end{equation}
has the eigenvalues $\big\{\xi_{j,k-1} \big\}_{j=1}^p$ and the eigenpairs $\left\{ {\left( {{\lambda _i},{x_i}} \right)} \right\}_{i = p + 1}^{2n}$. For the $k$-step, we need to find the feedback vectors $f_k, \ g_k \in \mathbb{R}^n$, by single-input feedback control
$$
M \ddot x(t) + C_k\dot x(t) + K_kx(t) = b_k\left(f_k^\top\dot x(t)e^{-\lambda \tau} + g_k^\top x(t)e^{-\lambda \tau}\right),
$$
such that the eigenvalues $\big\{\xi_{j,k-1} \big\}_{j=1}^p$ of the system \eqref{eq3.20} are reassigned to $\big\{\xi_{j,k} \big\}_{j=1}^p$ and the remaining eigenstructure keeps the no spill-over property. After $m$-step assignment, the feedback matrices $F$ and $G$ can be obtained. Then finding the feedback vectors $f_k, \ g_k \in \mathbb{R}^n$ can be expressed as the following problem.

\vspace{0.2cm}
\textbf{Problem 2.} \emph{Given $M$, $C$, $K$, $B$, $\Lambda_{1}$, $X_{1}$, $\tau$ and $\Sigma_{1}$. For $k = 1, 2,\ldots, m$, let $\{\xi_{jk}\}_{j=1}^p$ and $\left\{C_k,K_k \right\}$ be defined by \eqref{eq3.18} and \eqref{eq3.19} respectively.  Find the two feedback vectors $f_{k}$ and $g_{k}$ in sequence such that the single-input closed-loop delayed pencil
$$
P_{\tau k}(\lambda) = \lambda^2M + \lambda\left (C_{k} -  b_kf_k^\top e^{-\lambda \tau}\right) + K_{k} -  b_kg_k^\top e^{-\lambda \tau}
$$
has the desired eigenvalues $\big\{\xi_{jk} \big\}_{j=1}^p$ and the eigenpairs $\left\{(\lambda_i, x_i) \right\}_{i=p+1}^{2n}$.}

 \begin{lemma}\label{lemma3.1}\cite{10} For $k = 1, 2,\ldots, m$, and arbitrary $\beta_{k} \in \mathbb{C}^p$, define
\begin{equation}\label{eq3.21}
f_k = MX_1\beta_k, \ g_k = \left(MX_1\Lambda_1 + CX_1 \right )\beta_k,
\end{equation}
then the single-input closed-loop delayed pencil $P_{\tau k}(\lambda)$ has the desired eigenpairs $\left\{(\lambda_i, x_i) \right\}_{i=p+1}^{2n}$, i.e., the no spill-over property is preserved.
\end{lemma}

Let $\big\{(\xi_{jk},y_{jk})\big\}_{j=1}^p$ be the eigenpairs of the $P_{\tau k}(\lambda)$, and define
$$
D_k = {\rm{diag}}\big(\xi_{1k},\xi_{2k}, \ldots ,\xi_{pk}\big),\ Y_k = \big[y_{1k},y_{2k}, \ldots ,y_{pk}\big],
$$
\begin{equation}\label{eq3.22}
H_k =\left(h_{ls}^{(k)}\right)= X_1^\top MY_kD_k + \Lambda_1X_1^\top MY_k + X_1^\top CY_k.
\end{equation}
From Lemma \ref{lemma3.1}, solving \textbf{Problem 2} is transformed into determining vector $\beta_{k}$ such that
\begin{equation}\label{eq3.99}
M{Y_k}{D_k^2} + C_{k}{Y_k}{D_k}- b_{k}f_{k}^\top{Y_k}{D_k}e^{-\tau D_{k} } + K_{k}{Y_k} - bg_{k}^\top{Y_k}e^{-\tau D_{k} } = 0.
\end{equation}
Substituting  $f_{k}, g_{k}$ defined in \eqref{eq3.21} into \eqref{eq3.99}, and according to the definition of $H_{k}$ in \eqref{eq3.22}, we have
\begin{equation}\label{eq3.23}
\begin{aligned}
M{Y_k}{D_k^2} + C_{k}{Y_k}{D_k} + K_{k}{Y_k} &= b_{k}\beta_{k}^\top(X_1^\top MY_kD_k + \Lambda_1X_1^\top MY_k + X_1^\top CY_k)e^{-\tau D_{k}  }\\
& = b_{k}   \beta_{k}^\top H_{k} e^{-\tau D_{k}}:= b_{k}\gamma_{k},
\end{aligned}
\end{equation}
where $\gamma_{k}  = \beta_{k}^\top H_{k} e^{-\tau D_{k}}$. Similarly, choose $\gamma_{k}=[1,1,\ldots,1] \in \mathbb{R}^{p}$, then
\begin{equation}\label{eq3.24}
\beta_{k}^\top H_{k} = [1,1,\ldots,1]e^{\tau D_{k} }.
\end{equation}
From \eqref{eq3.23}, we can solve the linear systems
\begin{equation}\label{eq3.25}
\left(\xi_{jk}^2M + \xi_{jk}C_k + K_k\right)y_{jk} = b_k, \ j = 1,2, \ldots ,p
\end{equation}
to get $Y_{k}$ by choosing the appropriate parameters $\{\eta_{jk}\}$ such that the coefficient matrix of linear systems \eqref{eq3.25} are nonsingular. Then the matrix $H_{k}$ can be formed from \eqref{eq3.22}, and $\beta_{k}$ in \eqref{eq3.21} can be obtained by solving the linear system \eqref{eq3.24}.

 Clearly, the method mentioned earlier necessitates the solution of at least $((m-1)\cdot p)$ $n$-order linear systems. When dealing with a second-order system with a large degree of freedom, the computational expenses can be quite substantial. However, in practical engineering scenarios, the requirement to reassign poles is usually limited, with $p\ll n$ as a general rule. To minimize computational costs, we establish a new connection among system matrices, undesired eigenvalues and their corresponding eigenvectors, desired eigenvalues, and control vectors. Subsequently, we develop an efficient algorithm for forming matrix $H_{k}$ without the need to calculate $Y_{k}$.

% On the other hand, in practical engineering problems, the number of poles that need to be reassigned is very small, in general, $p\ll n$. In order to reduce the computational costs, we establish the new connection among system matrices, undesired eigenvalues and corresponding eigenvectors, desired eigenvalues and control vectors, and develop an efficient algorithm for forming matrix $H_{k}$ without calculating $Y_{k}$.

Based on the definition of $H_{k}$ in \eqref{eq3.22}, we can derive
%\begin{small}
\begin{displaymath}
\begin{aligned}
h_{ls}^{(k)} & = x_{l}^\top\left[\xi_{sk}M+\lambda_{l}M+C\right]y_{sk}     \\
             & = \frac{x_l^\top}{\xi_{sk}-\lambda_l}\left[\xi_{sk}^2M+\xi_{sk}C-\Big(\lambda_l^2M+\lambda_lC\Big)\right]y_{sk}  \\
             & = \frac{x_l^\top}{\xi_{sk}-\lambda_l}\Bigg(\xi_{sk}^2M+\xi_{sk}\Big(C-\sum\limits_{i=1}^{k-1} b_if_i^\top e^{-\tau D_{k}  }\Big)
                 +\Big(K-\sum\limits_{i=1}^{k-1} b_ig_i^\top e^{-\tau D_{k}  }\Big)            \\
             &\ \ \ \ \ \ -\Big(\lambda_l^2M+\lambda_lC+K\Big)  +\xi_{sk}\sum\limits_{i=1}^{k-1}b_if_i^\top e^{-\tau D_{k}  }+\sum\limits_{i=1}^{k-1}b_ig_i^\top e^{-\tau D_{k}  } \Bigg)y_{sk}  \\
             & =\frac{x_l^\top}{\xi_{sk}-\lambda_l}\Bigg(\Big(\xi_{sk}^2M+\xi_{sk}C_k+K_k\Big)-\Big(\lambda_l^2M+\lambda_lC + K\Big)
                +\sum\limits_{i=1}^{k-1}b_i\Big(\xi_{sk}f_i^\top+g_i^\top\Big)e^{-\tau D_{k}  }\Bigg)y_{sk}.
\end{aligned}
\end{displaymath}
%\end{small}
According to $\left(\lambda_l^2M+\lambda_lC + K\right)x_l=0$ and \eqref{eq3.25}, we have
\begin{equation}\label{eq3.26}
h_{ls}^{(k)} = \frac{x_l^\top b_k}{\xi_{sk}-\lambda_l} + \frac{e^{-\tau \xi_{sk}  }x_l^\top\left(\sum\limits_{i=1}^{k-1} b_i\left(\xi_{sk}f_i^\top+g_i^\top\right) \right)y_{sk}}{\xi_{sk}-\lambda_l}.
\end{equation}
Let
$$
G_k = X_1^\top\sum\limits_{i=1}^{k-1} b_i\beta_i^\top, \ W_k =\left(w_{ls}^{(k)}\right)= X_1^\top\sum\limits_{i=1}^{k-1} b_i\beta_i^\top H_k=G_kH_k.
$$
Substituting $H_{k}$ define in \eqref{eq3.22} into $W_{k}$, we can obtain
\begin{equation}\label{eq3.27}
\begin{aligned}
w_{ls}^{(k)} & = x_l^\top\Bigg(\sum\limits_{i=1}^{k-1} {b_i\beta_i^\top}\Big(\xi_{sk}X_1^\top M + \Lambda_1X_1^\top M + X_1^\top C\Big) \Bigg)y_{sk} \\
             & = x_l^\top\Bigg(\sum\limits_{i=1}^{k-1} b_i\Big(\xi_{sk}f_i^\top + g_i^\top\Big) \Bigg)y_{sk}.
\end{aligned}
\end{equation}
 \begin{remark}\label{remark3}From \eqref{eq3.26} and \eqref{eq3.27}, we find that each element of $W_{k}$ has a strong correlation with the corresponding element of $H_{k}$. We can use this relationship to derive a low order matrix equation about $H_{k}$.
 \end{remark}
 Let
$$
V_k = \left( {\begin{array}{*{20}{c}}
\frac{1}{\xi_{1k} - \lambda_1} & \ldots & \frac{1}{\xi_{pk} - \lambda_1} \\
            \vdots             & \ddots & \vdots                         \\
\frac{1}{\xi_{1k} - \lambda_p} & \ldots & \frac{1}{\xi_{pk} - \lambda_p}
\end{array}} \right),\
$$
$$
A_k={\rm{diag}}\left(x_1^\top b_k, \ldots, x_p^\top b_k\right), \
T_{k} = {\rm {diag}}\left(e^{-\tau \xi_{1k}}, \ldots,e^{-\tau \xi_{pk}}\right),\
R_{k} = V_{k}T_{k}.
$$
By \eqref{eq3.26} and \eqref{eq3.27}, $H_{k}$ satisfies the matrix equation described by
\begin{equation}\label{eq3.28}
H_k = U_k + (V_k T_{k}) \ast (G_kH_k)=U_k + R_k \ast W_k,
\end{equation}
where $U_k =\left(u_{ls}^{(k)}\right)= \left(\frac{x_l^Tb_k}{\xi_{sk}-\lambda_l} \right)=A_k V_k$ and $\ast$
represents the Hadamard product.

 Let $A= \big[a_1, a_2, \ldots, a_n\big ]$ be an $m \times n$ matrix, $a_{i}$ denote the $i$-th column of the matrix $A$, the $mn$ dimension vector ${\rm{vec}}(A)=\Big[a_1^\top$, $a_2^\top, \ldots, a_n^\top \Big]^\top$ be called the column straightening of matrix $A$, $A \otimes B$ represent the Kronecker product of matrix $A$ and $B$, and $I_{n}$ represent the $n$-order identity matrix. Then the column straightening of the matrix with the Hadamard product and Kronecker product of the matrix has the following properties.

 \begin{lemma}\label{lemma3.2}\cite{17} Let $A$ be an $m \times n$ matrix.\\
(1) If $B$ is an $m \times n$ matrix, then ${\rm{vec}}\left(A \ast B\right)={\rm {diag}}\left({\rm {vec}}(A)\right){\rm {vec}} (B)$;\\
(2) If $B$ and $C$ are $n \times p$ and $p \times q$ matrix, respectively, then ${\rm{vec}}(ABC)=(C^\top \otimes A){\rm {vec}}(B)$.
\end{lemma}

 From Lemma \ref{lemma3.2}, the matrix equation \eqref{eq3.28} can be transformed into the following $p^{2}$-order linear system
\begin{equation}\label{eq3.29}
\big[I_{p^2}-{\rm {diag}}\left({\rm{vec}}(R_k)\right)(I_p \otimes G_k)\big]{\rm{vec}}(H_k)={\rm{vec}}(U_k).
\end{equation}

 Let $Z_{k}=I_{p^2}-{\rm {diag}}\left({\rm{vec}}(R_k)\right)(I_p \otimes G_k)$, because $R_k$, $G_k$ and $U_k$ are all known,  the matrix $H_{k}$ can be obtained by solving the $p^2$-order linear system \eqref{eq3.29}, and  the $\beta_{k}$ can be obtained by solving the $p$-order linear system \eqref{eq3.24}. In summary, we can derive the following theorem.

\begin{theorem}\label{th2}For $k=1,2,\ldots m$, if the single-input control system $\left( M, C_{k},K_{k} , b_{k} \right)$ is partially controllable with respect to the $\big\{\xi_{j,k} \big\}_{j=1}^p$, $\left\{ {{\lambda _i}} \right\}_{i = 1}^p$ and $\big\{\xi_{jk} \big\}_{j=1}^p$ are all distinct, $\big\{\xi_{jk} \big\}_{j=1}^p \cap \left\{\lambda_i \right\}_{i=1}^{2n} = \emptyset, \
\left\{\lambda_i \right\}_{i = 1}^p \cap \left\{\lambda_i \right\}_{i=p+1}^{2n} = \emptyset$, and $\big\{y_{jk} \big\}_{j=1}^p$ are the solutions of \eqref{eq3.25}, then \textbf{Problem 2} has a solution, which can be expressed by \eqref{eq3.21}, where $\beta_{k}$ is the solution of linear system \eqref{eq3.24}.
\end{theorem}

 Based on the Theorem \ref{th1} and Theorem \ref{th2}, the calculation process of solving the multi-input partial quadratic pole assignment problem with time delay can be summarized as the following algorithm.

 \vspace{0.2cm}
\begin{tabular}{l}
\hline
{\bf Algorithm 2} A multi-step method for \textbf{Problem 1} by multi-input control.\\
\hline
{\bf Input: } \\
\indent  The symmetric matrices $M,C\in \mathbb{R}^{n\times n}$ with $M$ positive definite;\\
\indent  The full column rank control matrix $B \in {\mathbb{R}^{n\times m}}$ and the time delay $\tau > 0$;\\
\indent  The self-conjugate eigenpairs $\left\{(\lambda_i,x_i)\right\}_{i=1}^p$ and a self-conjugate set $\{ {{\mu_i}} \}_{i = 1}^p$.\\

{\bf Output: } \\
\indent State feedback real matrices $F$ and $G$. \\
1: Form $\Lambda_1 = {\rm{diag}}\big(\lambda_1,\lambda_2, \ldots ,\lambda_p\big)$ and $X_1 = \big[x_1,x_2, \ldots ,x_p\big]$;        \\
2: Choose $\eta_{jk}$ and compute $\xi_{jk} = \eta_{jk} + \frac{k}{m}\big(\mu_j - \eta_{jk}\big)$, $j = 1,2, \ldots ,p$, $k = 1,2, \ldots ,m$; \\
3: Compute $f_1 = MX_1\beta_1$, $g_1 = \big(MX_1\Lambda_1 + CX_1\big)\beta_1$ by Algorithm 1; \\
4: For $k = 2,3,\ldots , m$, do             \\
\indent 4.1: Compute $R_k,G_k,U_k$. If ${\rm{cond}}\left(Z_{k}\right)$ is large, take another $\eta_{jk}$;            \\
\indent 4.2: Compute $H_k$ by solving the linear system \eqref{eq3.29};    \\
\indent 4.3: Compute $\beta_k$ by solving the linear system \eqref{eq3.24}; \\
\indent 4.4: Compute $f_k = MX_1\beta_k$, $g_k = \left(MX_1\Lambda_1 + CX_1\right)\beta_k$; \\
5: Form $F = \big[f_1,f_2, \ldots ,f_m\big]$, $G = \big[g_1,g_2, \ldots ,g_m\big]$.\\
\hline
\end{tabular}

\section{ Numerical Examples }\label{example}

In this section, we provide three examples to illustrate the effectiveness of Algorithm 1 and Algorithm 2. All the numerical examples are implemented in MATLAB R2020a and run on the PC with 3.80 GHz Inter Core i7, 10700K CPU. We use $${\rm {Error_1}}=\big\|MY_1\Sigma_1^2 + CY_1\Sigma_1 - BF^\top Y_1\Sigma_1{\rm{e}}^{-\tau\Sigma_{1} }  + KY_1    - BG^\top Y_1{\rm{e}}^{-\tau\Sigma_{1} } \big\|_F$$
and
$${\rm {Error_2}}=\big\|M{X_2}\Lambda_2^2 + C{X_2}{\Lambda_2}- BF^\top{X_2}{\Lambda_2}{\rm{e}}^{-\tau\Lambda_{2} } + K{X_2} - BG^\top {X_2}{\rm{e}}^{-\tau\Lambda_{2} }\big\|_F$$
 to measure the residuals of the closed-loop system. Example 1 and Example 2 show the correctness of  Algorithm 2, and Example 3 shows Algorithm 2 is more efficient than the traditional multi-step method in \cite{15} in the case $p\ll n$.

  \vspace{0.2cm}
{\bf{Example 1.}} \cite{15} Consider the second-order control system
$$M=I_{3},\ C = \left({\begin{array}{*{20}{c}}
2.5  & 2 & 0 \\
2 & 1.7  &0.4 \\
  0  & 0.4   & 2.5
\end{array}} \right),\
K = \left( {\begin{array}{*{20}{c}}
16  & 12 &  0   \\
12  & 13 & 4   \\
 0  & 4 &  29
\end{array}} \right).
$$
Time delay $\tau = 0.1$ and $B = b = \big( {\begin{array}{*{20}{c}}
1 & 3 & 3
\end{array}} \big)^\top$. The corresponding open-loop pencil has 6 eigenvalues: $\lambda_{1,2} = - 0.0129 \pm 1.4389i$, $\lambda_{3,4} = - 1.3342 \pm 5.2311i$, $\lambda_{5,6} = - 2.0030 \pm 4.7437i$.
In order to make the system more stable, the first two eigenvalues $\lambda_{1,2} = - 0.0129 \pm 1.4389i$ are reassigned to $\mu_{1} =- 0.2$ and $\mu_{2} =- 0.3$, and the remaining eigenvalues of the original system are kept unchanged.  Using Algorithm 1, we can obtain $\beta_1 = ( - 0.3023 + 0.4678i , \ - 0.3023 - 0.4678i)^\top,$  the feedback vectors
$$F = f = \left( {\begin{array}{*{20}{c}}
0.1428  \\
-0.1541 \\
0.0215
\end{array}} \right),\
G = g = \left( {\begin{array}{*{20}{c}}
-0.9698  \\
1.2224 \\
-0.1852
\end{array}} \right),
$$
and
$$
{\rm {Error_1}} = 6.0497{\rm{e}}-15,\ {\rm {Error_2}} = 1.9486{\rm{e}}-13.
$$
 \vspace{0.2cm}

\begin{figure}[htbp]
\centering
\includegraphics[width=0.7\textwidth]{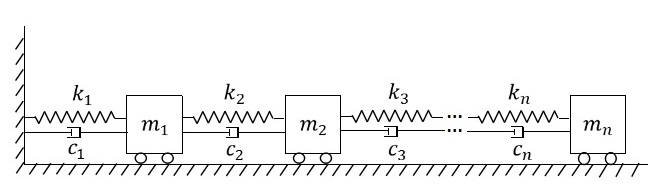}
\caption{Series mass-damping-spring system.}\label{fig1}
\label{fig:digit}
\end{figure}

{\bf{Example 2.}} Consider the same system matrices in Example 1 with time delay $\tau = 0.1$ and control matrix $B = \left( {\begin{array}{*{20}{c}}
1&2  \\
3&2\\
3&4
\end{array}} \right).$
The first two eigenvalues $\lambda_{1,2} =-0.0129 \pm 1.4389i$ are reassigned to $\mu_{1} =-0.2$ and $\mu_{2} =- 0.3$, and the remaining eigenvalues of the original system
are kept unchanged. Let $\eta_{jk} = \lambda_{j}, j=1,2,\ldots,p$. Applying Algorithm 2, we can obtain $\beta_1 = ( - 0.3779 + 0.4117i ,\  - 0.3779 - 0.4117i)^\top$, $\beta_2 = ( - 0.4108 - 0.3047i , \ - 0.4108 + 0.3047i)^\top$, the feedback matrices
$$F = \left( {\begin{array}{*{20}{c}}
{0.0220}&{ - 0.6561}\\
{ - 0.0131}&{ 0.7658}\\
{ 0.0005}&{ - 0.1141}
\end{array}} \right), \
G = \left( {\begin{array}{*{20}{c}}
{ - 1.0119}&{ - 0.2284}\\
{1.2347}&{ 0.0669}\\
{-0.1844}&{ 0.0047}
\end{array}} \right),
$$
and
$$
{\rm {Error_1}} = 1.5638{\rm{e}}-12, \  {\rm {Error_2}} = 2.0668{\rm{e}}-13.
$$
 \vspace{0.2cm}
\textbf{Example 3.} Consider a mass-damp-spring system with $n$ degrees of freedom shown as in Figure \ref{fig1}. The system matrices are $M ={\rm{diag}}\big(m_1,m_2,\ldots,m_n\big),$
\\
$$C = \left( {\begin{array}{*{20}{c}}
{{c _1} + {c _2}}&{ - {c _2}}&{}&{}\\
{ - {c _2}}&{{c _2} + {c _3}}&{ - {c _3}}&{}\\
{}& \ddots & \ddots &{}\\
{}&{}&{ - {c _n}}&{{c _n}}
\end{array}} \right), \
K = \left( {\begin{array}{*{20}{c}}
{{k _1} + {k _2}}&{ - {k _2}}&{}&{}\\
{ - {k _2}}&{{k _2} + {k _3}}&{ - {k _3}}&{}\\
{}& \ddots & \ddots &{}\\
{}&{}&{ - {k _n}}&{{k _n}}
\end{array}} \right).
$$
\indent Let $M=I_{n}, c_{1}=k_{1}=0, c_{i}=8, k_{i}=150, i = 2,3, \ldots ,n$, time delay $\tau = 0.1$ and control matrix $B = \big[e_1,e_2\big]$, where $e_1$ and  $e_2$ represent the first and second columns of matrix $I_n$, respectively. Let the degrees of freedom of the system be 500, 1000, 2000, 3000, 4000 and 5000. In each of the six cases, only one eigenvalue is unstable, and they are $\lambda_{969}=7.8049{\rm{e}}-14$, $\lambda_{1993}=6.4293{\rm{e}}-13$, $\lambda_{3997}=1.3911{\rm{e}}-12$, $\lambda_{5999}=1.5965e-12$, $\lambda_{7999}=1.9989{\rm{e}}-07$ and $\lambda_{9999}=1.6973{\rm{e}}-07$, respectively.
These eigenvalues are all reassigned to $\mu_1 = -0.2$. Let $\eta_{jk} = \lambda_{j}, j=1,2,\ldots,p$, the numerical results are shown in Table \ref{table1}.
\begin{table}
\begin{center}
\begin{tabular}{c|cccccc}
\hline
\multirow{2}{*}{${\rm{n}}$}&\multicolumn{2}{c}{${\rm{CPU \ Time(s)}}$}&\multicolumn{2}{c}{${\rm{Error_1}}$}&\multicolumn{2}{c}{${\rm{Error_2}}$} \\
\cline{2-7}
&Alg. 2&Alg. \cite{15}&Alg. 2&Alg. \cite{15}&Alg. 2&Alg. \cite{15}\\
\hline
500&0.0018&0.0095&1.8677e-13&1.8947e-13&4.5030e-09&4.5030e-09\\
\hline
1000&0.0029&0.0324&3.2900e-12&1.8293e-13&6.3745e-09&6.3745e-09\\
\hline
2000&0.0081&0.1677&2.5098e-12&1.9489e-13&1.0236e-08&1.0236e-08\\
\hline
3000&0.0252&0.3583&5.7325e-12&1.8936e-13&1.2699e-08&1.2699e-08\\
\hline
4000&0.0291&0.6484&2.5455e-12&1.9244e-13&1.4229e-08&1.4229e-08\\
\hline
5000&0.0415&1.0438&3.2718e-12&1.8946e-13&1.6358e-08&1.6358e-08\\
\hline
\end{tabular}
\end{center}
\caption{Example 3 numerical results}\label{table1}
\end{table}

Algorithm 2 shows that the calculations related to the degree of freedom $n$ are all the multiplication of matrices and vectors, and do not involve solving $n$-order linear systems. Specifically, we only need to solve the $p$-order linear system \eqref{eq3.24} and the $p^2$-order linear system \eqref{eq3.29}. Therefore, the time consumption of Algorithm 2 is much less than that of Algorithm in \cite{15}, and the results in Table \ref{table1} are consistent with the theoretical prediction.  Moreover,  as the degree of freedom  of the system increases, the time consumption of Algorithm \cite{15} increases markedly, while the time consumption of Algorithm 2 does not increase much.

\section{ Conclusion }\label{con}
In this paper, the partial pole assignment problem in symmetric quadratic pencil with time delay is considered. We  present a simple explicit solution for this problem for the single-input system, which involves only the multiplication of matrices and vectors. For the multi-input system, by establishing a new matrix equality relation, we develop a novel multi-step method to transform this problem into solving linear systems with low order.  Numerical examples show that the proposed method is effective, and both the computational time and cost of the method proposed in this paper are  markedly reduced than the traditional multi-step method  in the case $p\ll n$. Our method is  particularly well-suited for the large-scale system where only a small part of the spectrum needs to be reassigned.

\section*{Declaration of Competing Interest}
The authors declare that they have no known competing financial interests or personal relationships that could have appeared to
influence the work reported in this paper.

\section*{Data availability}
Data will be made available on request.

\section*{Acknowledgments}
The author sincerely thanks  Professor  Hua Dai for the careful guidance on  the partial pole assignment in symmetric quadratic pencil.

\section*{References}
\balance

\end{document}